\documentclass[10pt]{elsart}

\journal{}

\usepackage[a4paper, left=0.8in,right=0.8in,top=0.8in,bottom=0.8in]{geometry}
\usepackage{amsmath}
\usepackage{amssymb}
\usepackage{leftidx}
\usepackage{url}

\DeclareMathOperator\sgn{sgn}

\begin{document}

\newtheorem {Theorem}{Theorem}
\newtheorem {Corollary}{Corollary}

\newcommand {\dd}{\textup{d}}
\newcommand{\RR}{\mathbb{R}}
\newcommand {\sts}[2]{\genfrac\{\}{0pt}{}{#1}{#2}}

\newcommand*\pFqskip{8mu}
\catcode`,\active
\newcommand*\W{\begingroup
        \catcode`\,\active
        \def ,{\mskip\pFqskip\relax}%
        \dopFq
}
\catcode`\,12
\def\dopFq#1#2#3{%
        W\left(\genfrac..{0pt}{}{#1}{#2};#3\right)%
        \endgroup
}

\begin{frontmatter}

 \title{On the generalization of the Lambert $W$ function}
 \author[mezo]{Istv\'an Mez\H{o}}
 \address[mezo]{Department of Mathematics,\\Nanjing University of Information Science and Technology,\\No.219 Ningliu Rd, Nanjing, Jiangsu, P. R. China}
 \ead{istvanmezo81@gmail.com}
 \thanks[mezo]{The research of Istv\'an Mez\H{o} was supported by the Scientific Research Foundation of Nanjing University of Information Science \& Technology, and The Startup Foundation for Introducing Talent of NUIST. Project no.: S8113062001}
 \author[baricz]{\'Arp\'ad Baricz}
 \address[baricz]{Department of Economics, Babe\c{s}-Bolyai University, 400591 Cluj-Napoca, Romania\\
Institute of Applied Mathematics, \'Obuda University, 1034 Budapest, Hungary}
 \ead{bariczocsi@yahoo.com}

\begin{abstract}The Lambert $W$ function, giving the solutions of a simple transcendental equation, has become a famous function and arises in many applications in combinatorics, physics, or population dyamics just to mention a few. In this paper we construct and study in great detail a generalization of the Lambert $W$ which involves some special polynomials and even combinatorial aspects.
\end{abstract}

\begin{keyword}
\MSC Lambert $W$ function, Stirling numbers, Laguerre polynomials, transcendence
\end{keyword}
\end{frontmatter}

\section{Introduction}

\subsection{The definition of the Lambert $W$ function}

The solutions of the transcendental equation
\begin{equation}
xe^x=a\label{defw}
\end{equation}
were studied by Euler and by Lambert before \cite{W}. The inverse of the function on the left hand side is called the Lambert function and is denoted by $W$. Hence the solution is given by $W(a)$. If $-\frac1e<a<0$, there are two real solutions, and thus two real branches of $W$ \cite{Veberic}. If we enable complex values of $a$, we get many solutions, and $W$ has infinitely many complex branches \cite{W,CJ,CJK,JHC}. These questions are discussed by Corless et al. in details \cite{W}. The question why do we use the letter ``W'' is discussed by Hayes \cite{Hayes} who dedicated a whole article to the question; see also \cite{CJK}.

The Lambert function appears in many physical and mathematical problems, see for example the survey of Corless et al. \cite{W}. Here we concentrate on the mathematical aspects, and we plan to write a survey on the physical applications in the near future.

A simple mathematical application connects $W$ to the distribution of primes via the Prime Number Theorem \cite{Visser}. A more involved problem -- which was our original motivation to generalize the Lambert $W$ -- is the asymptotic estimation of the Bell and generalized Bell numbers. We discuss this question in more detail.

\subsection{A problem from combinatorics}

The $n$th Bell number $B_n$ counts in how many ways can we partition $n$ e\-le\-ments into nonempty subsets. Moser and Wyman \cite{MW} proved that in the asymptotic approximation of the Bell numbers the $W$ function appears. A simple expression due to Lov\'asz \cite[Section 1.14, Problem 9.]{Lovasz} states that,
\[B_n\sim \frac{1}{\sqrt{n}} \left( \frac{n}{W(n)} \right)^{n + \frac{1}{2}} \exp\left(\frac{n}{W(n)} - n - 1\right).\]
See also the note of Canfield \cite{Canfield}.

There is a generalization of the Bell numbers, called $r$-Bell numbers. The $n$th $r$-Bell number $B_{n,r}$ gives that in how many ways can we partition $n$ elements into nonempty subsets such that the first $r$ elements go to separated blocks \cite{MezoB}. The asymptotic behavior of the $r$-Bell numbers as $n\to\infty$ and $r$ is fixed was described by C. B. Corcino and R. B. Corcino \cite{CC}. They proved that in the asymptotic expression for $B_{n,r}$ one needs to solve the equation
\begin{equation}
xe^x+rx=n.\label{rLamdefeq}
\end{equation}

We generalize and rewrite this equation in the following steps:
\begin{equation}
xe^x+rx=n\longrightarrow (x-t)e^{cx}+r(x-t)=a\label{trans}
\end{equation}
\[\longrightarrow e^{-cx}=\frac{x-t}{a-r(x-t)}\longrightarrow e^{-cx}=-\frac1r\frac{x-t}{x-\left(\frac ar+t\right)}.\]
Hence, if we introduce the new variables
\[a_0=-\frac 1r,\quad\mbox{and}\quad s=\frac ar+t,\]
we arrive at the transcendental equation
\begin{equation}
e^{-cx}=a_0\frac{x-t}{x-s}.\label{geneq3}
\end{equation}

The above example shows that it is necessary to study the extension of \eqref{defw}, and \eqref{geneq3} shows a possible way to do it.

\section{The extensions of the equation defining $W$}

Beyond the above mathematical example the Scott et al. mentioned \cite{SBDM} that in some recent physical investigations on Bose-Fermi mixtures the transcendental equation
\begin{equation}
e^{-cx}=a_0\frac{(x-t_1)(x-t_2)\cdots(x-t_n)}{(x-s_1)(x-s_2)\cdots(x-s_m)}\label{defmw}
\end{equation}
appears. Motivated by these reasons, we initiate an investigation of these equations in general.

It is obvious that the parameter $c$ (if it is not zero) can be built in the parameters $t_i$ and $s_i$, by substituting $x/c$ in place of $x$:
\[e^{cx}\frac{(x-t_1)(x-t_2)\cdots(x-t_n)}{(x-s_1)(x-s_2)\cdots(x-s_m)}\longrightarrow c^{m-n}e^x\frac{(x-ct_1)(x-ct_2)\cdots(x-ct_n)}{(x-cs_1)(x-cs_2)\cdots(x-cs_m)}\]
Hence giving the solution of \eqref{defmw} is nothing else but evaluating the (generally multivalued) \emph{inverse} of the function
\[e^x\frac{(x-t_1)(x-t_2)\cdots(x-t_n)}{(x-s_1)(x-s_2)\cdots(x-s_m)}.\]

The first step could be a notation. We shall denote the inverse of this function in the point $a$ by
\[\W{t_1,t_2,\dots,t_n}{s_1,s_2,\dots,s_m}{a}.\]
We know that it is close to the hypergeometric function notation, but up to our present knowledge even the usual Lambert $W$ does not have anything to do with the hypergeometric functions.

So, in particular,
\[\W{}{}{a}=\log(a),\quad\W{0}{}{a}=W(a),\]

Apart from the trivial simplification
\[\W{t,t_1,t_2,\dots,t_n}{t,s_1,s_2,\dots,s_m}{a}=\W{t_1,t_2,\dots,t_n}{s_1,s_2,\dots,s_m}{a}\]
there is no obvious reducibility of this general $W$ function. In what follows we will examine two special cases,
\[\quad\W{t}{s}{a},\quad\mbox{and}\quad\W{t_1,t_2}{}{a},\]
We restrict us to the real line, however, complex extension of this function class would be interesting.

In the followings we study the properties of some particular generalizations. Physical applications of this generalizations can be found in \cite{MezoKeady}

\section{The case of one upper and one lower parameter}

\subsection{Basic properties}

First we look for the Taylor series of $\W{t}{s}{a}$ around $a=0$. This question reveals a surprising connection between this function and the Laguerre polynomials.

\begin{thm}\label{TaylorWts}The Taylor series of $\W{t}{s}{a}$ around $a=0$ is
\[\W{t}{s}{a}=t-T\sum_{n=1}^\infty\frac{L_n'(nT)}{n}e^{-nt}a^n,\]
where $T=t-s\neq 0$, and $L_n'$ is the derivative of the $n$th Laguerre polynomial \cite{Bell,Szego}.
\end{thm}

\textit{Proof.} We use the Lagrange Inversion Theorem \cite[p. 14.]{AS} -- a primary tool in inverse function investigations. Let
\[f(x)=e^x\frac{x-t}{x-s}.\]
We choose a point in which $f$ is zero, and then we invert its series in this point. This point is $t$. Hence, by Lagrange's theorem,
\begin{equation}
\W{t}{s}{a}=t+\sum_{n=1}^\infty\frac{a^n}{n!}\lim_{w\to t}\frac{\dd^{n-1}}{\dd w^{n-1}}\left(\frac{w-t}{f(w)}\right)^n.\label{WtsLag}
\end{equation}
Here the only difficulty is the expression
\begin{equation}
\lim_{w\to t}\frac{\dd^{n-1}}{\dd w^{n-1}}\left(\frac{w-t}{f(w)}\right)=\lim_{w\to t}\frac{\dd^{n-1}}{\dd w^{n-1}}\left(\frac{w-s}{e^w}\right).\label{tmp1}
\end{equation}
One thing we can forecast to make the things easier. Recalling the Rodrigues formula of the generalized Laguerre polynomials \cite{AS,Szego}
\[L_n^{(\alpha)}(x)=\frac{x^{-\alpha}e^x}{n!}\frac{\dd^n}{\dd x^n}(e^{-x}x^{n+\alpha}),\]
it can be seen that a modification of \eqref{tmp1} will lead to a generalized Laguerre polynomial. Let us make this precise. The limit \eqref{tmp1}, if we expand it entirely, takes the form
\begin{equation}
e^{-nt}\sum_{k=1}^nA_{n,k}(-1)^{k-1}T^k,\label{exprT}
\end{equation}

and the coefficients $A_{n,k}$ we determine now. An easy calculation gives the following table:
\[
\begin{array}{c|ccccccc}A_{n,k}&k=1&k=2&k=3&k=4&k=5&k=6&k=7\\\hline
n=1 & 1 & \text{} & \text{} & \text{} & \text{} & \text{} & \text{}\\
n=2 & 2 & 2 & \text{} & \text{} & \text{} & \text{} & \text{} \\
n=3 & 6 & 18 & 9 & \text{} & \text{} & \text{} & \text{}\\
n=4 & 24 & 144 & 192 & 64 & \text{} & \text{} & \text{}\\
n=5 &120 & 1200 & 3000 & 2500 & 625 & \text{} & \text{}\\
n=6 &720 & 10800 & 43200 & 64800 & 38880 & 7776 & \text{}\\
n=7 &5040 & 105840 & 617400 & 1440600 & 1512630 & 705894 & 117649\\
\end{array}
\]
These number directly cannot be found in the Online Encyclopedia of Integer Sequences (OEIS) \cite{OEIS} -- a prominent tool of investigators facing to some unknown integer sequence -- but the row sum appears under the identification number A052885. There we can find that $A052885$ equals to the sum of
\[(n-1)!n^{k-1}\binom{n}{k}\frac{1}{(k-1)!}.\]
Hence we can suspect that this is exactly what we are looking for,
\[A_{n,k}=(n-1)!n^{k-1}\binom{n}{k}\frac{1}{(k-1)!}.\]
Once we have this conjecture, the proof is easy (by induction). Hence we can step forward, substituting this into \eqref{WtsLag}:
\[\W{t}{s}{a}=t+\sum_{n=1}^\infty\frac{(ae^{-t})^n}{n}\sum_{k=1}^n\frac{(-n)^{k-1}}{(k-1)!}\binom{n}{k}T^k=\]
\[\W{t}{s}{a}=t-\sum_{n=1}^\infty\frac{(ae^{-t})^n}{n^2}\sum_{k=1}^n\binom{n}{k}\frac{(-nT)^k}{(k-1)!}.\]
Recalling the explicit expression for the generalized Laguerre polynomials \cite[p. 775, Table 22.3]{AS}:
\[L_n^{(\alpha)}=\sum_{k=0}^n\binom{n+\alpha}{n-k}\frac{(-x)^k}{k!},\]
we can easily see that our inner sum in the Taylor series is simply
\[\sum_{k=1}^n\binom{n}{k}\frac{(-nT)^k}{(k-1)!}=-nTL_{n-1}^{(1)}(nT).\]
The relation \cite[p. 778]{AS}
\[L_{n-1}^{(1)}(x)=-L_n'(x),\]
finalizes the proof. Note that $T\neq0$ is necessary in the above argument, especially in \eqref{exprT}. $T=0$ stands for the standard $\log$ function.
\hfill\qed

We could not find the radius of convergence of the above Taylor series in general, just when $T<0$.

\begin{thm}\label{TaylorWtsconvrad}When $t<s$ the radius of convergence of the power series in Theorem \ref{TaylorWts} is the following
$$r=e^{t}\lim_{n\to\infty}\left|\frac {L_{n-1}^{(1)}(nT)} {L_n^{(1)}((n+1)T)}\right|=e^{\frac{t+s}{2}-2\sqrt{s-t}}.$$
\end{thm}

\textit{Proof.} By using the asymptotic behavior of Laguerre polynomials for large $n$, but fixed $\alpha$ and $x>0$, that is,
$$L_{n}^{\alpha}(-x)=\frac{(n+1)^{\frac{\alpha}{2}-\frac{1}{4}}}{2\sqrt{\pi}}\frac{e^{-\frac{x}{2}}}{x^{\frac{\alpha}{2}+\frac{1}{4}}}e^{2\sqrt{x(n+1)}}\left(1+\mathcal{O}\left(\frac{1}{\sqrt{n+1}}\right)\right),$$
we have that if $T<0,$ that is, $t<s,$ then
$$\frac {L_{n-1}^{(1)}(nT)} {L_n^{(1)}((n+1)T)}=e^{-\frac{T}{2}-2\sqrt{-T}}\left(\frac{n}{n+1}\right)^{-\frac{1}{2}}\frac{1+\mathcal{O}\left(\frac{1}{\sqrt{n}}\right)}{1+\mathcal{O}\left(\frac{1}{\sqrt{n+1}}\right)},$$
which implies that the radius of convergence in this case is $r=e^{\frac{t+s}{2}-2\sqrt{s-t}}.$
\hfill\qed

\subsection{A straight generalization of the classical $W$ function -- The $r$-Lambert function}

The definition of $\W{t}{s}{a}$ shows that this function is not a direct generalization of the Lambert $W$ function in the sense that there are no finite special values of $t$ and $s$ for which $\W{t}{s}{a}$ would be equal to $W(a)$. On the other hand, we saw that there is another equation, namely \eqref{rLamdefeq}, which solution easily specializes to $W$ and, via the transformation \eqref{trans}, to $\W{t}{s}{a}$ as well. Hence equation \eqref{rLamdefeq} deserves a deeper look. And there is one more reason: the function $e^x\frac{x-t}{x-s}$ has singularity while $xe^x+rx$ has not, so the analysis of this latter function is easier.

Let $r$ be a fixed real number. Motivated by the work of Corcino et al. \cite{CC}, the inverse of the function $xe^x+rx$ we call $r$-Lambert function and we denote it by $W_r$. The steps under \eqref{trans} yield the following observations.

\begin{thm}The solution(s) of the equation
\[e^{cx}\frac{x-t}{x-s}=a\]
is (are)
\[x=t+\frac1cW_{-ae^{-ct}}\left(cae^{-ct}T\right).\]
Hence, in particular,
\[\W{t}{s}{a}=t+W_{-ae^{-t}}\left(ae^{-t}T\right).\]
Here $T=t-s$.

Moreover, the equation -- which is a generalization of \eqref{rLamdefeq} and so of \eqref{defw} --
\[(x-t)e^{cx}+r(x-t)=a\]
can be resolved by the $r$-Lambert function as
\[x=t+\frac1cW_{re^{-ct}}\left(cae^{-ct}\right).\]
\end{thm}

Depending on the parameter $r$, the $r$-Lambert function has one, two or three real branches and so the above equations can have one, two or three solutions (we restrict our investigation to the real line). Now we describe this in more detail when $r\neq0$.

\begin{thm}\label{Branchprops}Depending on the value of $r$, we can classify $W_r(x)$ as follows.
\begin{enumerate}
	\item If $r>1/e^2$, then $W_r(x):\RR\to\RR$ is a strictly increasing, everywhere differentiable function such that $\sgn(W_r(x))=\sgn(x)$.

	\item If $r=1/e^2$, then $W_r(x):\RR\to\RR$ is a strictly increasing function which is differentiable on $\RR\setminus\{-4/e^2\}$. Moreover, $\sgn(W_r(x))=\sgn(x)$.
	\item If $0<r<1/e^2$, then there are three branches of $W_r$ which we denote by $W_{r,-2}$, $W_{r,-1}$, $W_{r,0}$. Let
	\[\alpha_r=W_{-1}(-re)-1,\quad\mbox{and}\quad\beta_r=W_0(-re)-1,\]
	where $W_{-1}$ and $W_0$ are the two branches of the Lambert function.
	
	Then $\alpha_r$ and $\beta_r$ determine the above branches as follows:	
\begin{itemize}
	\item $W_{r,-2}:\,]-\infty,f_r(\alpha_r)]\to\,]-\infty,\alpha_r]$ is a strictly increasing function,
	\item $W_{r,-1}:\,[f_r(\alpha_r),f_r(\beta_r)]\to\,[\alpha_r,\beta_r]$ is a strictly decreasing function, and
	\item $W_{r,0}:\,[f_r(\beta_r),+\infty[\,\to\,[\beta_r,+\infty[$ is a strictly increasing function.
\end{itemize}
These three functions are differentiable on the interior of their domains.

	\item Finally, if $r<0$, then $W_r$ has two branches, $W_{r,-1}$ and $W_{r,0}$. Let
	\[\gamma_r=W_(-re)-1,\]
	where $W$ is the classical Lambert function.
	
	Then for these branches we have that
\begin{itemize}
	\item $W_{r,-1}:[f_r(\gamma_r),+\infty[\,\to]-\infty,\gamma_r]$ is strictly decreasing, while
	\item $W_{r,0}:[f_r(\gamma_r),+\infty[\,\to[\gamma_r,+\infty[$ is a strictly increasing function.
\end{itemize}
Both of them are differentiable in the interior of their domains.
	
\end{enumerate}
Here $\sgn$ is the signum function such that $sgn(0)=0$; and $f_r(x)=xe^x+rx$.
\end{thm}

The branches $W_{r,-2}$ and $W_{r,-1}$ take negative values, while $W_{r,0}$ is positive whatever values have $r$ provided that it is smaller than $1/e^2$. Hence the notation.

\textit{Proof.} The proof is entirely elementary. Which we have to take into account is that the derivative of $f_r(x)=xe^x+rx$ equals to zero exactly in $W(-re)-1$:
\[f_r'(x)=e^x(1+x)+r=0\Longrightarrow e^{1+x}(1+x)=-re\Longrightarrow x=W(-re)-1.\]
Then analyzing the possible values of this $x$ we can easily confirm the number, domain and range of the different branches. The monotonicity and differentiability follows from the similar properties of $f_r=W_r^{-1}$. We shall see immediately why the point $W_{1/e^2}(-2)$ is exceptional.
\hfill\qed

Now we look for the derivative and integral of $W_r$.

\begin{thm}The derivative of the $r$-Lambert function is
\begin{equation}
W_r'(x)=\frac{1}{e^{W_r(x)}(1+W_r(x))+r}.\label{Wdiff}
\end{equation}
In particular, taking into account that $W_r(0)=0$,
\[W_r'(0)=\frac{1}{1+r}.\]
The integral of $W_r(x)$ is
\[\int W_r(x)\dd x=\frac{r}{2}W_r^2(x)+e^{W_r(x)}(1-W_r(x)+W_r^2(x))+c.\]
\end{thm}

\textit{Proof}. We can mimic the proof of Corless et. al \cite{W}.
\hfill\qed

Since $f_r(-1)=-\frac1e-r$ and $f_r(1)=e+r$, we have the special values
\[W_r\left(-\frac1e-r\right)=-1,\quad\mbox{and}\quad W_r(e+r)=1.\]
(One can also see that the equation
\[W(x\log x)=\log x\]
generalizes to
\[W_r[(x+r)\log x]=\log x\quad(x>0).)\]

By the above theorem we have that
\[\int_0^{r+e}W_r(x)\dd x=\frac{r}{2}+e-1.\]

Also, this theorem helps us to find the point where $W_r$ is not differentiable. The denominator in \eqref{Wdiff} is zero in
\[W(-re)-1.\]
In these points we always have branch cut except when $r=1/e^2$ where the two branches can be continuously connected together under the cost that $W_{1/e^2}$ is not differentiable at
\[f_{1/e^2}(W(-1/e^2\cdot e)-1)=f_{1/e^2}(-2)=-\frac4{e^2}.\]
However, at this point $W_{1/e^2}$ is continuous (as everywhere on the real line) and
\[W_{1/e^2}\left(-\frac4{e^2}\right)=-2.\]

In the following part we reveal a very interesting relationship between the $r$-Lambert function and combinatorics.

\subsection{The $r$-Lambert function and the Stirling numbers}

Around $x=0$ the Lambert function has the Taylor series expansion
\begin{equation}
W(x)=\sum_{n=1}^\infty\frac{(-n)^{n-1}}{n!}x^n\label{Wseries}
\end{equation}
with convergence radius $\rho=1/e$. The sequence $n^{n-1}$ is the number of rooted trees on $n$ labelled points. Hence $W$ has combinatorial connections (not just by this reason, see \cite[(4.18)]{W} or \cite{GKP,JCHK} for connections to Stirling numbers of the first kind).

Now we point out that the series \eqref{Wseries} nicely generalizes not just keeping its combinatorial meaning but extending it. To do this, we need some basic facts from finite set partition theory. The symbol $\sts{n}{k}$ denotes the Stirling number of the second kind with parameters $n$ and $k$. The number $\sts{n}{k}$ gives that in how many ways can be partitioned a set of $n$ elements into $k$ nonempty, disjoint subsets. The exponential generating function for a fixed $k$ is
\begin{equation}
\sum_{n=k}^\infty\sts{n}{k}\frac{x^n}{n!}=\frac{1}{k!}(e^x-1)^k.\label{genst}
\end{equation}
The following polynomial identity we also will need:
\begin{equation}
x^n=\sum_{k=1}^n\sts{n}{k}x^{\underline k},\label{polidst2}
\end{equation}
where
\[x^{\underline k}=x(x-1)\cdots(x-k+1)\quad (x^{\underline 0}=1)\]
is the falling factorial. We shall use the rising factorial, too:
\[x^{\overline k}=x(x+1)\cdots(x+k-1)\quad (x^{\overline 0}=1).\]
This latter is often denoted by $(x)_k$ and called as Pochhammer symbol. We adopt the notations of \cite{GKP} where the reader can find more information on Stirling numbers.

First we show that the Taylor series of the $r$-Lambert function contains these numbers and then we reveal connections of this series to another combinatorial problem.

\begin{thm}\label{TaylorrLambert}Let $r\neq-1$ be a fixed real number. Moreover, we define the polynomial $M_k^{(n)}(y)$ as
\begin{equation}
M_k^{(n)}(y)=\sum_{i=1}^kn^{\overline i}\sts{k}{i}(-y)^i.\label{Mkdef}
\end{equation}
Then in a neighbourhood of $x=0$,
\begin{equation}
W_r(x)=\frac{x}{r+1}+\sum_{n=2}^\infty M_{n-1}^{(n)}\left(\frac{1}{r+1}\right)\frac{x^n}{(r+1)^nn!}.\label{TaylorWr}
\end{equation}
\end{thm}

Firstly, it is important to note that when $r=-1$, by the fourth point of Theorem \ref{Branchprops}, $\gamma_{-1}=W_0(e)-1=0$ is a branch point, so the Taylor series does not exist. (Analytic continuation would resolve the problem.)

Secondly, let us realize that if $r=0$, then
\[M_{n-1}^{(n)}(1)=\sum_{i=1}^{n-1}n^{\overline i}\sts{n-1}{i}(-1)^i=\sum_{i=1}^{n-1}(-n)^{\underline i}\sts{n}{i}=(-n)^{n-1}\]
by \eqref{polidst2}. Hence we get back the well known Taylor expansion of the Lambert function.

\textit{Proof.} By the Lagrange inversion we can write $W_r(x)$ around $x=0$ as
\begin{equation}
W_r(x)=\sum_{n=1}^\infty\frac{x^n}{n!}\lim_{w\to0}\frac{\dd^{n-1}}{\dd w^{n-1}}\left(\frac{1}{e^w+r}\right)^n\label{WrTaylorpr}
\end{equation}
on some neighbour of the origin. Although the $(n-1)$-th derivative would be sufficient, we determine
\[\frac{\dd^{k}}{\dd w^{k}}\left(\frac{1}{e^w+r}\right)^n\]
for every $1\le k\le n$. This will lead to nice combinatorial results. We can rewrite this $k$th derivative as
\begin{equation}
\frac{\dd^{k}}{\dd w^{k}}\left(\frac{1}{e^w+r}\right)^n=\frac{1}{(e^w+r)^n}\sum_{i=1}^kC_{i,k}^{(n)}\frac{e^{iw}}{(e^w+r)^i}.\label{Cikdiff}
\end{equation}
Letting $w\to0$
\[W_r(x)=\sum_{n=1}^\infty\frac{x^n}{(r+1)^nn!}\sum_{i=1}^{n-1}\frac{C_{i,n-1}^{(n)}}{(r+1)^i}.\]
Hence we do not need to do more just find the coefficients ${C_{i,k}^{(n)}}$ where $i=1,\dots,k$ and $k=1,\dots,n$. First, we can easily see that
\[{C_{1,1}^{(n)}}=-n,\quad\mbox{and that}\quad{C_{1,k+1}^{(n)}}={C_{1,k}^{(n)}}.\]
Or, together,
\begin{equation}
{C_{1,k}^{(n)}}=-n\quad(k=1,2,\dots).\label{invalcik}
\end{equation}
Therefore if we write ${C_{i,k}^{(n)}}$ in a triangle with a fixed $n$, indexing the columns by $i$ and the rows by $k$, we have that the first column is $-n$. By induction, the general recurrence in this triangle reads as:
\begin{equation}
C_{i,k+1}^{(n)}=iC_{i,k}^{(n)}-(n+i-1)C_{i-1,k}^{(n)}\quad(i=1,2,\dots,k+1).\label{recurcik}
\end{equation}
This can be directly seen comparing
\[\frac{\dd^{k}}{\dd w^{k}}\left(\frac{1}{e^w+r}\right)^n\quad\mbox{with}\quad\frac{\dd^{k+1}}{\dd w^{k+1}}\left(\frac{1}{e^w+r}\right)^n.\]
Using the initial values \eqref{invalcik}, this recurrence gives all the coefficients in principle. However, at this point we still cannot realize the real combinatorial meaning of these coefficients. To get a deeper insight, we determine the exponential generating function
\begin{equation}
A_i^{(n)}(x):=\sum_{k=i}^\infty C_{i,k}^{(n)}\frac{x^k}{k!}.\label{Aidef}
\end{equation}
Our recurrence \eqref{recurcik}  yields that $A_i^{(n)}(x)$ must satisfy the differential equation
\[\frac{\dd}{\dd x}A_i^{(n)}(x)=iA_i^{(n)}(x)+(n+i-1)A_{i-1}^{(n)}(x).\]
The first function is
\[A_1^{(n)}(x)=\sum_{k=1}^\infty C_{1,k}^{(n)}\frac{x^k}{k!}=-n\sum_{k=1}^\infty\frac{x^k}{k!}=-n(e^x-1)\]
by \eqref{invalcik}. Now the differential equation for $A_2^{(n)}(x)$ says that
\[\frac{\dd}{\dd x}A_2^{(n)}(x)=2A_2^{(n)}(x)-(n+1)A_1^{(n)}(x).\]
Substituting $A_1^{(n)}(x)$,
\[\frac{\dd}{\dd x}A_2^{(n)}(x)=2A_2^{(n)}(x)+n(n+1)(e^x-1).\]
This, together with the initial value $A_2^{(n)}(0)=0$ (see \eqref{Aidef}) results that
\begin{align*}
A_2^{(1)}(x)&=1-2e^x+e^{2x},\\
A_2^{(2)}(x)&=3-6e^x+3e^{2x},\\
A_2^{(3)}(x)&=6-12e^x+6e^{2x},
\end{align*}
and, in general,
\[A_2^{(n)}(x)=\frac12n(n+1)(e^x-1)^2\quad(n=1,2,\dots).\]
Similarly,
\begin{align*}
A_3^{(n)}(x)&=-\frac16n(n+1)(n+2)(e^x-1)^3\quad(n=1,2,\dots),\\
A_4^{(n)}(x)&=\frac1{24}n(n+1)(n+2)(n+3)(e^x-1)^4\quad(n=1,2,\dots).
\end{align*}
The pattern is obvious:
\[A_i^{(n)}(x)=\frac{(-1)^i}{i!}n^{\overline i}(e^x-1)^i\quad(i,n=1,2,\dots).\]
Applying the generating function \eqref{genst} and comparing the coefficients, we arrive at the nice and simple expression for $C_{i,k}^{(n)}$:
\[C_{i,k}^{(n)}=(-1)^in^{\overline i}\sts{k}{i}.\]
Based on \eqref{Cikdiff}, we now have that
\[\lim_{w\to0}\frac{\dd^{n-1}}{\dd w^{n-1}}\left(\frac{1}{e^w+r}\right)^n=\frac{1}{(r+1)^n}\sum_{i=1}^kn^{\overline i}\sts{k}{i}\left(\frac{-1}{r+1}\right)^i.\]
This together with \eqref{WrTaylorpr} provides our result.
\hfill\qed

\subsubsection{The polynomial $M_k^{(n)}(y)$}

It is worth to study the polynomials $M_k^{(n)}(y)$ a bit more profoundly. The generating function of these polynomials reveals that this polynomial is connected not just to the Stirling numbers but to another combinatorial sequences as well.
\begin{equation}
\sum_{k=1}^\infty M_k^{(n)}(y)\frac{x^k}{k!}=\sum_{k=1}^\infty\left(\sum_{i=1}^k C_{i,k}^{(n)}y^i\right)\frac{x^k}{k!}\label{egenfMk}
\end{equation}
\[=\sum_{i=1}^\infty y^i\sum_{k=i}^\infty C_{i,k}^{(n)}\frac{x^k}{k!}=\sum_{i=1}^\infty A_i^{(n)}(x)y^i\]
\[=\sum_{i=1}^\infty\frac{(-1)^i}{i!}n^{\overline i}(e^x-1)^iy^i=-1+\frac{1}{(1+(e^x-1)y)^n}.\]
Thus, in particular, $y=1$ (that is, $r=0$) gives that
\[\sum_{k=1}^\infty M_k^{(n)}(1)\frac{x^k}{k!}=-1+e^{-nx},\]
and so
\[M_k^{(n)}(1)=(-n)^k\quad(k>0).\]
Setting $k=n-1$ as it is in Theorem \ref{TaylorrLambert}, we get back the classical Taylor series of the Lambert function again.

If we set $y=-1$ ($r=-2$), then
\begin{equation}
\sum_{k=1}^\infty M_k^{(n)}(-1)\frac{x^k}{k!}=-1+\frac{1}{(2-e^x)^n}.\label{Mkgenf}
\end{equation}
For $n=1$ this series has well known coefficients in combinatorics:
\[\frac{1}{2-e^x}=\sum_{n=1}^\infty F_n\frac{x^n}{n!},\]
where $F_n$ is the $n$th ordered Bell number, preferential arrangement, or Fubini number \cite{James,Sprugnoli,Tanny,VC}. Hence
\[M_k^{(1)}(-1)=F_k\quad(k\ge1).\]
(Note that $M_0^{(1)}(-1)=0$ while $F_0=1$ by convention.) Of course, this is also obvious from the \eqref{Mkdef} definition of $M_k^{(n)}(y)$ and from the fact that
\[F_n=\sum_{k=1}^nk!\sts{n}{k}.\]
A somewhat greater surprise that $M_k^{(n)}(-1)$ counts, for a general $n$, the barred preferential arrangements recently studied by Ahlbach, Usatine and Pippenger.

A preferential arrangement on $n$ elements is a partition of these elements into groups (blocks) such that the order of the groups counts but the elements in the individual blocks does not count. The number of all the possible preferential arrangements on $n$ elements is exactly the $n$th Fubini number $F_n$. In addition, a barred preferential arrangement with parameter $k$ is a preferential arrangement such that $k$ bars are placed to separate the blocks into $k+1$ sections. The total number of the barred preferential arrangements on $n$ elements with a fixed $k$ is given by the numbers $r_{k,n}$ studied by the above mentioned three authors in detail \cite{AUP}.

Comparing the generating functions in \cite[Theorem 4.]{AUP} and \eqref{Mkgenf}, the values of $M_k^{(n)}(-1)$ coincide with these numbers:
\[M_k^{(n)}(-1)=r_{n-1,k}\quad(n,k\ge1).\]

\subsection{Generalizing the Omega constant}

The omega constant is the solution of the transcendental equation
\[\Omega e^\Omega=1.\]
Equivalently,
\[e^{-\Omega}=\Omega,\quad\mbox{or}\quad\log(\Omega)=-\Omega.\]
It is also the value of the Lambert function at 1:
\[\Omega=W(1).\]
Having the $r$-Lambert function, one can consider such constants by considering the transcendental equations determining $W_r$. Maybe the simplest constant after $\Omega$ one can define is the constant given by the equation
\[\Omega_1e^{\Omega_1}+\Omega_1=1,\]
that is, the value of $W_1$ at 1:
\[\Omega_1=W_1(1).\]
(More precisely, this function is $W_{1,0}$, the rightmost branch.)

The numerical value of $\Omega_1$ is
\[\Omega_1\approx0.401058137541547\dots\]

This number is transcendental. Let us suppose, in contrary, that $\Omega_1$ is al\-geb\-raic. By definition,
\[e^{\Omega_1}=\frac{1}{\Omega_1}-1,\]
hence $e^{\Omega_1}$ is algebraic. This is a contradiction by the Lindemann-Weierstrass theorem \cite{Baker}.

We note that the Lindemann-Weierstrass theorem gives us much more about the transcendence of different values of $W_r$.

\begin{thm}If $a$ and $r$ are algebraic numbers such that $a\neq0$, then $W_r(a)$ is transcendental.
\end{thm}

The Taylor series \eqref{TaylorWr} -- which still converges at $x=1$ -- gives that this constant can be represented as a rapidly convergent infinite sum:
\[\Omega_1=\frac12+\sum_{n=2}^\infty\frac{1}{2^nn!}\sum_{k=1}^{n-1}n^{\overline k}\sts{n-1}{k}\left(-\frac12\right)^k.\]

After these explorations we turn to the convergence questions of the series \eqref{TaylorWr}.

\subsection{About the radius of convergence for the Taylor series of $W_r$}

\subsubsection{The $r=-2$ case}

In \cite{AUP} one can find the asymptotic estimation for $r_{k,n}$, which helps us to find the radius of convergence of the Taylor series of $W_{-2,-1}(x)$ (see the next point why the Taylor series belongs to this branch indexed by $-2$). Namely, a remark after Theorem 14. in \cite{AUP} says that
\[r_{k,n}\sim\frac{(n+k)!}{2^{k+1}k!(\log(2))^{k+n+1}},\]
from which an estimation comes for $M_{n-1}^{(n)}(-1)=r_{n-1,n-1}$. By using Stirling's formula, the radius of convergence of the Taylor series \eqref{TaylorWr} can be determined:
\begin{equation}
\rho_{-2}=\frac{\log^2(2)}{2}\approx0.24023\label{rhom2}
\end{equation}

\subsubsection{The $r<0$ case in general}

The saddle point method \cite{Wilf} will help us to say something about $\rho_r$ in general. The exponential generating function \eqref{egenfMk} of $M_k^{(n)}(y)$ has (real) singularity in the following cases:
\begin{enumerate}
	\item $\frac{1}{r+1}=y<0$, that is, $r<-1$, or
	\item $\frac{1}{r+1}=y>1$, that is, $-1<r<0$.
\end{enumerate}

Then standard arguments from the saddle point method say that if $r<0$, then $M_k^{(n)}(y)$ (where, as always, $y=1/(r+1)$) surely grows faster than $k!a^n$ for some real $a$ ($k!$ comes from the fact that we are talking about an exponential generating function). In particular,
\[M_{n-1}^{(n)}(y)\gg(n-1)!a^n\quad(n\to\infty),\]
and so
\[M_{n-1}^{(n)}(y)\frac{1}{(r+1)^nn!}\gg\frac{(n-1)!a^n}{(r+1)^nn!}=\frac{\left(\frac{a}{r+1}\right)^n}{n}\quad(n\to\infty).\]
This shows that in this case \eqref{TaylorWr} has a finite convergence radius $\rho_r$. How large can be $\rho_r$? The exact answer is hard to find, because it is hard to estimate $M_{n-1}^{(n)}(y)$ (the more simpler Bell polynomials do not contain rising factorials like $M_{n-1}^{(n)}(y)$ and it is hard to estimate them).

If $r<-1$, then $\gamma_r>0$ and thus the \eqref{TaylorWr} Taylor series belongs to the left branch $W_{r,-1}$ while if $0>r>-1$, this series represents (a part of) $W_{r,0}$.

\subsubsection{The $r>0$ case}

We know that $\rho_0=\frac1e$. How $\rho_r$ varies when $r>0$? Since in this case the exponential generating function of $M_k^{(n)}(1/(r+1))$ still has singularity, $\rho_r$ must be finite.

We have to clarify to which branch does the Taylor series belong. There are three branches, $W_{r,-2}$, $W_{r,-1}$, and $W_{r,0}$. The branch cuts can be read out from the third point of Theorem \ref{Branchprops}: $f_r(\alpha_r)$ and $f_r(\beta_r)$. Since the domain of $W_{r,0}$ starts at $f_r(\beta_r)<0$, the Taylor series around 0 surely belongs to $W_{r,0}$.

\subsection{The asymptotics of the $r$-Lambert function at $\pm\infty$}

It is known \cite{W} that
\[W(x)\sim\log(x)+\log\left(\frac{1}{\log(x)}\right)\]
as $x$ grows. We shall prove the corresponding result for $W_r(x)$.

\begin{thm}\label{Wrapprox}For any $r\in\RR$
\[W_r(x)\sim\log(x)+\log\left(\frac{1}{\log(x)}-\frac{r}{x}\right)\]
as $x\to\infty$.

Moreover,
\[W_r(x)\sim\frac1rx\]
as $x\to-\infty$.
\end{thm}

\textit{Proof.} We follow the ideas described in \cite{W}. Let us rewrite $W_r(x)$ as
\begin{equation}
W_r(x)=\log(x)+u(x).\label{defu}
\end{equation}
By using the definition of $W_r(x)$,
\[W_r(x)e^{W_r(x)}+rW_r(x)=x,\]
we have
\[(\log(x)+u(x))xe^{u(x)}+r(\log(x)+u(x))=x.\]
Since $|u(x)|\ll\log(x)$ if $x$ is large, we have approximately that
\[\log(x)xe^{u(x)}+r\log(x)=x,\]
i.e.,
\[e^{u(x)}=\frac{\frac{x}{\log(x)}-r}{x}.\]
By taking logarithm, the first part of the theorem comes. The second part is easier, since $xe^x$ is negligable comparing to $rx$ if $x$ tends to $-\infty$, hence $xe^x+rx$ ``simplifies'' to $rx$. And this latter has the inverse $\frac1rx$. \hfill\qed

Note that the first asymptotic estimation of the previous theorem simplifies to the known case when $r=0$, but the second one does not have corresponding classical version. This is not a surprise, because the domain of the two real valued branches of the Lambert function do not extend to values $<-1/e$.

Note also that the above theorem helps to find exact asymptotics for the $r$-Bell numbers which was originally expressed in terms of the $r$-Lambert function in the work of Corcino et al. \cite{CC}.

We know that $W_r(x)\to0$ as $x\to0$ for any fixed $r$. The order of this convergence would be interesting. If $x$ is close to 0, we can use the Taylor series of Theorem \ref{TaylorrLambert}. The problem is that the approximation of the polynomials $M_k^{(n)}(y)$ is rather difficult.

\section{The case of two upper parameters}

The function $\W{t_1,t_2}{}{a}$ is the inverse of
\[e^x(x-t_1)(x-t_2).\]
Scott, Mann, and Martinez \cite{SMM} proposed a solution for the equation
\begin{equation}
a=e^{Rx}(x-t_1)(x-t_2).\label{SMM}
\end{equation}
as a product of two Lambert function of the form $cW(A)W(B)$, where $c$, $A$ and $B$ are determined by the parameters of the equation.

Before introducing our own solution, we point out that the Scott-Mann-Martinez (SMM) solution is not satisfactory. Indeed, depending on the values of $A$ and $B$, there are three possible cases: $W(A)W(B)$ can have one, two or four different real values. These happen when
\[A,B>0,\]
\[\quad-\frac1e<A<0, B>0,\quad-\frac1e<B<0, A>0,\]
or
\[-\frac1e<A,B<0,\]
respectively (considering the two real branches of $W$).

However, a short analyis of \eqref{SMM} shows that it can have one, two or three solutions. Hence for a class of parameters the SMM solution does not give back all the solutions and, in other cases, it gives a virtual solution. One might study exactly when the SMM solution works properly.

Similar steps as under \eqref{trans} show that \eqref{SMM} can be transformed to
\[c^2a=e^x(x-ct_1)(x-ct_2),\]
and so the solution of \eqref{SMM} is given by
\[\W{ct_1,ct_2}{}{c^2a}.\]

Since the work of Scott, Mann, and Martinez comes from physical problems to be solved (quantum physics and general relativity), it is worth to study the function $\W{t_1,t_2}{}{a}$ in its generality.

The Taylor series of $\W{t}{s}{a}$ (see Theorem \ref{TaylorWts}) contained the derivative of the Laguerre polynomials. In the case of $\W{t_1,t_2}{}{a}$ we have that the Bessel polynomials \cite{KF}
\[B_n(z)=\sum_{k=0}^n\frac{(n+k)!}{k!(n-k)!}\left(\frac{z}{2}\right)^k\]
appear (sometimes $y_n(z)$ notation is used). The following is a result of Mugnaini \cite{Mugnaini}.

\begin{thm}\label{TaylorWtt}The Taylor series of $\W{t_1,t_2}{}{a}$ around $a=0$ is
\begin{equation}
\W{t_1,t_2}{}{a}=t_1-\sum_{n=1}^\infty\frac{1}{n n!}\left(\frac{ane^{-t_1}}{T}\right)^nB_{n-1}\left(\frac{-2}{nT}\right),\label{Mugnainiresult}
\end{equation}
where $T=t_2-t_1$.
\end{thm}

Neither Mugnaini nor us could calculate the radius of convergence of this series.

\section{Closing remarks}

Closing the paper we list some additional questions which can serve and deserve further investigations.

\begin{itemize}
	\item The complex Lambert $W$ function has a nice but not elementary branch structure. It would be very interesting to find the branch structure of the functions introduced in the present paper, especially of the complex $r$-Lambert function $W_r$.

	\item We could not say nothing about the radius of convergence of \eqref{TaylorWr} for the $r$-Lambert function except the only case $r=-2$. What can we say about the convergence radius $\rho_r$ in general?
	
	\item There is a derivation formula for the composite function $W(e^x)$ in which the second order Eulerian numbers appear \cite[(3.5)]{W}. Does exist such formula for the $r$-Lambert function involving some generalization of the Eulerian numbers? If so, do these new numbers have some combinatorial meaning? See also \cite{Kruchinin}.
	
	\item After Theorem \ref{Wrapprox} we noted that we do not know how quickly $W_r(x)$ tends to zero. This can be a question to be answered in the future.
	
	\item Theorems \ref{TaylorWts} and \ref{TaylorWtt} contain Taylor series including the derivative of Laguerre polynomials, and an ``unidentified'' polynomial sequence. Could we say more on these Taylor series, especially find the radius of convergence in general apart from Theorem \ref{TaylorWtsconvrad}? Can we connect the $P_n$ polynomials to more known ones?
	
	\item There are good approximations to the Lambert function \cite{HH,Stewart}. How these generalize to the $r$-Lambert function?
	
	\item Could one carry out some general analysis for $\W{t_1,\dots,t_n}{s_1,\dots,s_m}{a}$?
	
	\item What is the radius of convergence of the Taylor series \eqref{Mugnainiresult}?
\end{itemize}

The first author wrote a C code to calculate the real $r$-Lambert function on all the branches. This can be downloaded from \url{https://sites.google.com/site/istvanmezo81/}.

\section{Acknowledgement}

The first author is very grateful to Prof. Grant Keady for his kind help during the preparation of the manuscript.


\begin{thebibliography}{99}
\bibitem{AS}
M. Abramowitz, I. A. Stegun, Handbook of Mathematical Functions with Formulas, Graphs, and Mathematical Tables, Dover, 1972.
\bibitem{AUP}
C. Ahlbach, J. Usatine, N. Pippenger, Barred preferential arrangements, Electron J. Combin 20(2) (2013) \#P55. Available online at

http://www.combinatorics.org/ojs/index.php/eljc/article/view/v20i2p55
\bibitem{Baker}
A. Baker, Transcendental Number Theory, Cambridge Univ. Press, 1990.
\bibitem{Bell}
W. W. Bell, Special Functions for Scientists and Engineers, D. Van Nostrand Company, 1968.
\bibitem{Canfield}
R. Canfield, The Moser-Wyman expansion of the Bell numbers, online note available at http://www.austinmohr.com/Work\_files/bellMoser.pdf
\bibitem{CC}
C. B. Corcino, R. B. Corcino, An asymptotic formula for $r$-Bell numbers with real arguments, ISRN Discrete Math, 2013 (2013), Article ID 274697, 7 pages.
\bibitem{W}
R. M. Corless, G. H. Gonnet, D. E. G. Hare, D. J. Jeffrey, D. E. Knuth, On the Lambert $W$ function, Adv. Comput. Math. 5 (1996), 329-359.
\bibitem{CJ}
R. M. Corless, D. J. Jeffrey, The unwinding number, \textsc{Sigsam} Bulletin 30(2) (1996), 28-35.
\bibitem{CJK}
R. M. Corless, D. J. Jeffrey, D. E. Knuth, A sequence of series for the Lambert $W$ function, ISSAC '97 Proceedings of the 1997 International Symposium on Symbolic and Algebraic Computation, 197-204.
\bibitem{GKP}
R. L. Graham, D. E. Knuth, O. Patashnik, Concrete Mathematics, Addison-Wesley, 1994.
\bibitem{Hayes}
B. Hayes, Why $W$?, American Scientist 93 (2005), 104-108.
\bibitem{HH}
A. Hoorfar, M. Hassani, Inequalities on the Lambert $W$ function and hyperpower function, J. Inequal. Pure and Appl. Math., 9(2) (2008), Art. 51. Available online at http://www.emis.de/journals/JIPAM/article983.html
\bibitem{James}
R. D. James, The factors of a square-free integer, Canad. Math. Bull. 11 (1968), 733-735.
\bibitem{JCHK}
D. J. Jeffrey, R. M. Corless, D. E. G. Hare, D. E. Knuth, Sur l'inversion de $y^\alpha e^y$ au moyen de nombres de Stirling asoci\'es, C. R. Acad. Sci. Paris, S\'erie I 320 (1995), 1449-1452.
\bibitem{JHC}
D. J. Jeffrey, D. E. G. Hare, R. M. Corless, Unwinding the branches of the Lambert $W$ function. Math. Scientist. 21 (1996), 1-7.
\bibitem{KF}
H. L. Krall, O. Frink, A new class of orthogonal polynomials: the Bessel polynomials, Trans. Amer. Math. Soc. 65(1) (1949), p. 100-115.
\bibitem{Kruchinin}
V. Kruchinin, Derivation of Bell polynomials of the second kind, ArXiv preprint, available at http://arxiv.org/pdf/1104.5065v1.pdf
\bibitem{Lovasz}
L. Lov\'asz, Combinatorial Problems and Exercises, North-Holland, 1993.
\bibitem{MezoB}
I. Mez\H{o}, The $r$-Bell numbers, J. Integer Seq. 14(1) (2011), Article 11.1.1.
\bibitem{MezoKeady}
I. Mez\H{o}, G. Keady, Some physical applications of generalized Lambert functions, ArXiv preprint, arxiv.org/abs/1505.01555
\bibitem{MW}
L. Moser, M. Wyman, An asymptotic formula for the Bell numbers, Trans. Royal Soc. Canada, Section III 49 (1955), 49-54.
\bibitem{Mugnaini}
G. Mugnaini, Generalization of Lambert $W$-function, Bessel polynomials and transcendental equations, ArXiv preprint, available at
http://arxiv.org/abs/1501.00138v2
\bibitem{SBDM}
T. C. Scott, J. F. Babb, A. Dalgarno, J. D. Morgan III, J. Chem. Phys. 99 (1993), 2841-2854.
\bibitem{SMM}
T. C. Scott, R. Mann, R. E. Martinez, General relativity and quantum mechanics: towards a generalization of the Lambert $W$ function, Appl. Algebra Engrg. Comm. Comput. 17(1) (2006), 41-47.
\bibitem{Stewart}
S. M. Stewart, On certain inequalities involving the Lambert $W$ function, J. Inequal. Pure and Appl. Math., 10(4) (2009), Art. 96. Available online at

\url{http://www.emis.de/journals/JIPAM/article1152.html}
\bibitem{OEIS}
The On-line Encyclopedia of Integer Sequences, http://oeis.org
\bibitem{Sprugnoli}
R. Sprugnoli, Riordan arrays and combinatorial sums, Discrete Math. 132 (1994), 267-290.
\bibitem{Szego}
G. Szeg\H{o}, Orthogonal Polynomials, fourth ed., Amer. Math. Soc, 1975.
\bibitem{Tanny}
S. M. Tanny, On some numbers related to the Bell numbers, Canad. Math. Bull. 17 (1975), 733-738.
\bibitem{Veberic}
D. Veberi\u{c}, Having fun with the Lambert $W(x)$ function, ArXiv preprint, available at http://arxiv.org/pdf/1003.1628.pdf
\bibitem{VC}
D. J. Velleman and G. S. Call, Permutations and combination locks, Math. Mag. 68(4) (1995), 243-253.
\bibitem{Visser}
M. Visser, Primes and the Lambert $W$ function, ArXiv preprint, available at http://arxiv.org/pdf/1311.2324v1.pdf
\bibitem{Wilf}
H. S. Wilf, Generatingfunctionology, Academic Press, 1994.
\end{thebibliography}
\end{document}